\documentclass[a4paper,twoside,12pt]{article}

\newtheorem{thm}{Theorem}

\usepackage[english]{babel}
\usepackage{babelbib}
\RequirePackage{geometry}
\usepackage[utf8]{inputenc}
\usepackage[T1]{fontenc}

\usepackage{comment}
\usepackage{float}
\usepackage{etoolbox}

\usepackage{txfonts} 

\usepackage[pdftex]{graphicx} 

\usepackage{listingsutf8}

\title{Automated conjecturing of Frobenius numbers via grammatical evolution}
\author{Nikola Adžaga \footnote{Faculty of Civil Engineering, University of Zagreb, e-mail: nadzaga@grad.hr}}
\date{}  

\usepackage{xcolor}
\usepackage{syntax} 

\def\<#1>{\synt{#1}}

\newcounter{grammarchild}
\newcommand{\grammartag}{\hfill(\refstepcounter{grammarchild}\textcolor{blue}{\thegrammarchild})}
\newenvironment{mygrammarrule}{\refstepcounter{grammarchild}\setcounter{grammarchild}{-1}\begin{grammar}}{\end{grammar}}
 
\usepackage{hyperref}
\hypersetup{
 unicode=true,
 colorlinks=true,
 linkcolor=black,
 citecolor=black,
 filecolor=black,
 urlcolor=black,
 bookmarksopen=true,
 bookmarksopenlevel=0,
 bookmarksnumbered=true,
 pdftitle={Automated conjecturing of Frobenius numbers via grammatical evolution},
 pdfauthor={Nikola Adžaga},
 pdfkeywords={},
 pdfstartpage={},
 pdfstartview=FitH,
 pdfnewwindow=true
}

\begin{document}
\maketitle

\section*{Abstract}
Conjecturing formulas and other symbolic relations occurs frequently in number theory and combinatorics. If we could automate conjecturing, we could benefit not only from speeding up, but also from finding conjectures previously out of our grasp. Grammatical evolution, a genetic programming technique, can be used for automated conjecturing in mathematics. Concretely, this work describes how one can interpret the Frobenius problem as a symbolic regression problem, and then apply grammatical evolution to it. In this manner, a few formulas for Frobenius numbers of specific quadruples were found automatically. The sketch of the proof for one conjectured formula, using lattice point enumeration method, is provided as well.

Same method can easily be used on other problems to speed up and enhance the research process.


\section{Introduction}
Unlike automation of theorem-proving, automated conjecturing is still a rare theme in mathematics. Probably the most significant example of a program used for automated conjecturing is Fajtlowicz's Graffiti \cite{graffiti}, which produced hundreds of conjectures in graph theory and inspired dozens of papers (e.\ g.\ \cite{Erdos}). Another example is Colton's HR \cite{Colton}. Colton considered mathematics as "a new domain for datamining" and used very simple algorithms on existing databases of mathematical knowledge, discovering, for example, a necessary condition for perfect numbers. Advances in artificial intelligence and its successes in other sciences suggest that, simply by treating examples and results of math experiments as data for machine learning problem, we could automatically gain new insights.

One method that could be very useful for mathematics is grammatical evolution (GE). GE is a genetic programming technique (i.\ e., used for automated programming), and as such, it can easily generate conjectures of deterministic type. On the other hand, it is a very general method, applicable to any problem which can be described using a context-free grammar. Indeed, it was successfully used in fields as various as computer-aided design (\cite{cad}) and computational finance (\cite{fin}). To author's knowledge, this is the first instance of using grammatical evolution in mathematics, so this work can be seen as a proof-of-concept demonstration -- grammatical evolution can be used for conjecturing in mathematics. We chose to show that for Frobenius problem. It is particularly suitable for this study because it consists of finding formulas, so GE is applied in its most usual form.

Frobenius problem is the following: given relatively prime natural numbers $a_1, \dots, a_n$, find the \emph{Frobenius number}, the largest natural number that is not representable as a non-negative integer linear combination of $a_1, \dots, a_n$. This number is denoted as $g(a_1, \dots, a_n)$.

For $n=2$, it is long known that $g(a_1, a_2) = a_1a_2-a_1-a_2$ (usually attributed to \cite{1882}).

For $n=3$, the quest for a similar simple formula has been the subject of extensive research. Curtis (\cite{msc12330}) has shown that this quest, in a way, cannot succeed -- there is no finite set of polynomials $\{ f_1, \dots, f_k \}$ such that, for each choice of $a_1, a_2, a_3$, there is some $i \in \{ 1, \dots, k \}$ such that $f_i(a_1, a_2, a_3) = g(a_1, a_2, a_3)$.
Still, efficient algorithms for calculating Frobenius number have been developed, as well as numerous lower and upper bounds (for an extensive review on these and other results regarding Frobenius number one can check \cite{alfonsinBook}).

For a fixed $n \geqslant 4$, there is a polynomial time algorithm, but it is impractical. There is no proven efficient algorithm, but various formulas for some specific cases were found, e.\ \nolinebreak g., for arithmetic and geometric sequences of arbitrary length. For $n=4$ more specific formulas were found, e.\ g.

\[ g(a, a+1, a+2, a+4) = (a+1)\left\lfloor \frac{a}{4} \right\rfloor + \left\lfloor \frac{a+1}{4} \right\rfloor + 2\left\lfloor \frac{a+2}{4} \right\rfloor - 1. \]

We show that grammatical evolution can be used for automated conjecturing in mathematics. The rest of this work is divided as follows. 
The second section is the most important part -- it describes how one can interpret the Frobenius problem as a symbolic regression problem, and then apply grammatical evolution to it. In this manner, a few formulas for Frobenius number were conjectured, including the following one:
\[ \displaystyle g(3k+1, 3k+4, 6k+3, 6k+9) = (3k+1)\left(k - \left\lfloor \frac{3k+1}{21} \right\rfloor \right) - 1 \quad \forall k \geqslant 5. \]
This section also shows the briefest sketch of the proof of this formula, using lattice point enumeration method. More important is the fact that, simply by changing the input data, we easily generated conjectures for quadruples of different kind, and even for sextuples.
The last section discusses the benefits and limitations of using grammatical evolution as automated conjecturing method in mathematics.



\section{The Frobenius problem as Symbolic Regression}
We interpret the Frobenius problem as a symbolic regression problem, and then approach this problem using grammatical evolution. Grammatical evolution is a method of automatic programming which uses grammars and evolutionary approach. Binary chromosomes are used. \nolinebreak{\emph{Codons}}, groups of $8$ bits, represent integer numbers. Codons determine which rule from the grammar will be used.

Following the short background on grammars, we explain how we used grammatical evolution and results we obtained.
\subsection{Context free grammars and Backus--Naur form (BNF)}
\label{BNF}
Grammars describing programming languages are defined following the 
formal languages, which were introduced as a scientific discipline by Noam Chomsky (e.\ g.~\cite{Chomsky1956}) 
during his search for simple rules capable of describing natural language.
Chomsky concluded that context-free grammars are sufficient for formalizing the grammar of the English language. These grammars form theoretical basis for most programming languages, but the importance of grammars is even greater since they provide the means for formalizing knowledge. If a finite set of simple rules can describe the grammar of English language, then such sets can describe many other phenomena, which in turn opens up a possibility for automated research of these phenomena. 


The Backus--Naur form is the most commonly used notation for expressing the grammar of a language. BNF defines syntax rules (as in~\cite{naur:algol60:cacm:1963}), 
and it formalizes syntactic expressions. By listing the set of production rules, we completely determine the grammar. In each rule there is exactly one variable on its left hand side, while on its right hand side there are expressions which can replace the given variable. 
The example of BNF is given in the next subsection (it is the grammar used in this research).

\subsection{Interpreting the Problem as Symbolic Regression Problem}
We calculated Frobenius number for $40$ quadruples of the form $(x, x+3, 2x+1, 2x+7)$ with a relatively simple recursive algorithm (using dynamic programming). These quadruples and their Frobenius numbers serve as the input data for grammatical evolution. Form of these quadruples was chosen as a simple form for which Frobenius number was previously unknown.

\newpage 
The input data then looks like this:
\begin{table}[H]
\begin{center}
{
     \begin{tabular}{ | c c c c | c |} 
	\hline
	x & x+3 & 2x+1 & 2x+7 & Frobenius number (target function)\\ \hline
	3 & 6 & 7 & 13 & 11 \\ \hline
	4 & 7 & 9 & 15 & 10 \\ \hline
	5 & 8 & 11 & 17 & 14  \\ \hline
	\dots & \dots & \dots & \dots & \dots \\ \hline
\end{tabular}}
\end{center}
\caption{Input data -- Frobenius number for quadruples $(x, x+3, 2x+1, 2x+7)$.}
\end{table}

A simple grammar was used (this is an example of BNF):
\begin{mygrammarrule}
<expr> ::= <expr> <op> <expr> \grammartag
\alt (<expr> <op> <expr>) \grammartag
\alt <expr> / <const> \grammartag
\alt (<expr> / <const>) \grammartag
\alt <var> \grammartag
\end{mygrammarrule}
\begin{mygrammarrule}
<op> ::= + \grammartag
\alt - \grammartag
\alt * \grammartag
\end{mygrammarrule}
\begin{mygrammarrule}
<var> ::= x \grammartag
\alt <const> \grammartag
\end{mygrammarrule}
\begin{mygrammarrule}
<const> ::= <const> <op> <const> \grammartag
\alt (<const> <op> <const>) \grammartag
\alt <const> / <const> \grammartag
\alt (<const> / <const>) \grammartag
\alt 1.0 \grammartag
\alt 3.0 \grammartag
\end{mygrammarrule}

Individual's chromosome is seen as a sequence of codons (integers). Example individual is $(120, 44, 42, 96, 189, 64, \ldots)$. We briefly explain how does this sequence map to an expression.  Start symbol is <expr>, the first symbol in the grammar. There are five rules given which can replace <expr>. First codon determines which rule is applied. Regardless of their integer value, codons are not skipped, but considered modulo the number of possibilities. First codon value is $120$, so rule number $0$ is applied, because $120 \equiv 0\ (\textrm{mod}\ 5)$. In this case, <expr> $\rightarrow$ <expr><op><expr>. Now, second codon determines which rule will be applied to <expr> (first of the three symbols). Second codon value is $44$, so the rule number $4$ is applied ($44 \equiv 4\ (\textrm{mod}\ 5)$), i.\ e., <expr><op><expr> $\rightarrow$ <var><op><expr>. Third codon value, $42$, determines that <var> will be replaced by $x$ (grammar gives just $2$ rules for <var> and $42 \equiv 0\ (\textrm{mod}\ 2)$, so we use the rule number $0$), i.\ e., <var><op><expr> $\rightarrow$ x<op><expr>. Now, $x$ is a terminal symbol (grammar does not provide any rules to replace it), so the next symbol we concentrate on is <op>. Next codon $96 \equiv 0\ (\textrm{mod}\ 3)$, so x<op><expr> $\rightarrow$ x+<expr>. The remaining two codons result in $x+\<expr>$ $\rightarrow$ $x+\<var>$ $\rightarrow$ $x+x$.

The aim is to find an individual that will, using this grammar, map to the expression describing the given data. A given sample of Frobenius numbers is used for the evaluation of indviduals -- fitness for this problem is given by the sum, taken over 40 function values, of the error between the evolved and target function  (Frobenius number). For more details on GE, one can check the first publication \cite{oneill:2001:TEC}. 

This search is done by evolutionary algorithm. As the individuals are simple bit-vectors, we do not have to employ any specific crossover or mutation operators. Standard genetic operators of one-point crossover and one-bit mutation were used. Parameters of the evolution were: population had $500$ individuals, it took $100$ generations, crossover probability was $0.9$, mutation probability was $0.1$. It is not important to carefully optimize these parameters.

GE did not find a formula valid for all given quadruples, but it found one which had a very high fit, i.\ e., a very small error:

$$g(x, x+3, 2x+1, 3x+7) \approx x(x/3 - x/21) - 1.$$
It was not hard to notice that this formula is valid in the case of $x = 3k+1$ and $k \geqslant 5$, which gives us the following result.
\begin{thm}
Frobenius number of the quadruple $(3k+1, 3k+4, 6k+3, 6k+9)$, for all $k \geqslant 5$, is
$$g(3k+1, 3k+4, 6k+3, 6k+9) = (3k+1)\left(k - \left\lfloor \frac{3k+1}{21} \right\rfloor \right) - 1. \quad $$ 
\end{thm}
(Floor function appears as integer division in C language which was used.)

\subsection{Proof of the conjectured formula}
The proof follows the method developed in the article \emph{Frobenius numbers by lattice point enumeration} \cite*{einstein2007frobenius} and concrete details were found with help of the algorithm, i.\ e., Mathematica package they developed. It is very similar to their proof for quadratic sequences. This is why author provides only the sketch ot the proof for the case of $k=7a$, without explaining the lattice point enumeration method (for understanding the proof and terminology, sections $1.-4., 6. \,\&\, 17.$ of \cite{einstein2007frobenius} should suffice).

Table $1$ is basically the proof for the formula when $k=7a$ ($a \in \mathbb{N}$). One can observe that the proposed \emph{elbows} determine the domain whose volume is $3k+1=21a+1$ (volume these elbows determine is a linear function of $a$, so it suffices to check two cases), while \emph{protoelbows} below show that the corresponding elbows are not in the fundamental domain.

 \begin{table}[H]
\begin{center}
{
     \begin{tabular}{ | p{2.8cm} | c c c |} 
     \hline
	Corner \mbox{weights} & $126a^2+27a$ & $105a+16$ & $105a+10$  \\ \hline 
	Corners & $(0, 0, 3a)$ & $(1, 1, 1)$ & $(1, 2, 0)$ \\ \hline 
	Elbows & $(0, 0, 3a+1)$ & $(0, 1, 2)$ $(2, 1, 0)$ & $(0, 3, 0)$ $(0, 2, 1)$ \\ \hline 
	Protoelbows & $(-2, 0, 3a+1)$ & $(-5, 1, 2)$ $(2, 1, -1)$ & $(-1, 3, 0)$ $(-3, 2, 1)$\\ \hline
     \end{tabular}

     \begin{tabular}{ | p{2.8cm} | c c c |} 
     \hline
	Corner \mbox{weights} &  $126a^2+27a-1$ & $126a^2+27a-2$  & $126a^2+27a-3$  \\ \hline
	Corners & $(2, 0, 3a-1)$ & $(4, 0, 3a-2)$ & $(6, 0, 3a-3)$  \\ \hline
	Elbows & $(1, 0, 3a)$ & $(3, 0, 3a-1)$ & $(5, 0, 3a-2)$ $\, (7, 0, 0)$ \\ \hline
	Protoelbows & $(1, -2, 3a)$ & $(3, -1, 3a-1)$ & $\,\,\,(0, 0, 0)$  $\quad (7, 0, -3)$ \\ \hline
	\end{tabular}
}
\end{center}
\caption{Description of fundamental domain for $(21a+1, 21a+4, 42a+3, 42a+9)$.}
\end{table}
\vskip 0.5em
The Frobenius corner here is $(0, 0, 3a)$, so Frobenius number is \[ 3a\cdot(6k+9)-(3k+1) = 3a\cdot(42a+9)-21a-1=126a^2+6a-1 = (21a+1)(7a-a)-1, \] \vskip -0.5em which is exactly what the formula gives for $k=7a$.

\subsection{Other conjectures}
Simply by changing the input data, system can generate other conjectures. 
One can try to generalize the formulas for arithmetic and geometric sequences. Probably the simplest generalization of these two is linear recursive sequence, e.\ g., sequence defined by $a_{n+1} = 3a_{n}+2$.
GE system resulted with some conjectures in these cases -- for odd numbers $x$ (the quadruples are not relatively prime for even numbers): 

\[ f(x = 4k+1, y = 3x+2, z = 3y+2, w = 3z+2) = \]
$$ 48k^2+16k-1-3(4k+1)\left(k+\left\lfloor\frac{2k}{13}\right\rfloor + \left\lfloor\frac{2k+6}{13}\right\rfloor - \left\lfloor\frac{2k+19}{26}\right\rfloor \right), $$

\[ f(x = 4k+3, y = 3x+2, z = 3y+2, w = 3z+2) = \]
$$ 48k^2+64k+19-3(4k+3)\left(k+\left\lfloor\frac{2k+1}{13}\right\rfloor + \left\lfloor\frac{2k+7}{13}\right\rfloor - \left\lfloor\frac{k+3}{13}\right\rfloor \right). $$

Even the size of the tuple is easily changed, e.\ g., system found the following conjecture as well:
\[ g(6k+1, 6k+4, 6k+7, 12k+3, 12k+9, 12k+15) = (6k+1)\left(k - \left\lfloor \frac{k}{13} \right\rfloor \right)-1 \quad \forall k \geqslant 5.  \]

One can even change the problem (again, simply by changing the input data to examples for that problem), but that would take us out of the scope of this work.

\section{Benefits and limitations}
Grammatical evolution can generate solutions in an arbitrary language by a simple change of the grammar -- which is usually one short input file. This advantage emerges from the generality of context-free grammars, hence GE has broad generativity -- it can be applied not only to the search of a formula, but also to anything context-free grammars can describe. 

Knowledge embedded in grammar's rules, speeds up the search process. In the present study, knowledge was basically the way expressions are created. Minor expectation of desired formula, the fact that division operator should appear only with constant as divisor, was also embedded in the grammar. The challenge in similar studies could be in embedding more complex expectations (e.\ g., recent results) in grammars. Constructing the grammar might require more creativity when one applies this method on different types of problems (say, constructing a grammar to find a graph with some particular property probably will not be such an easy task). However, given that usage of grammars is exactly the reason of great generalizability, the author believes that this advantage outweighs disadvantages.

Formulas found here could have been found using a more standard type of regression or using polynomial interpolation. However, because of the floor function, the formula proven here would be broken into seven polynomial formulas. Far more important, these methods require researcher to know or guess the model prior to applying them, and they are not so easy to generalize. On the other hand, fine tuning, i.\ e., getting the constants right, is a potential, but minor problem for GE. Formulas like the one proven here, ending with $-1$, can be missed by that $-1$ (in many runs of the program). Unlike computer, researcher can easily note and fix this type of problem.

\bibliographystyle{babplain} 

\bibliography{ge}

\end{document}